\documentclass[a4paper]{article}
\usepackage{amsfonts}
\usepackage{amssymb}
\usepackage{amsmath}
\usepackage{latexsym}
\newtheorem{theorem}{Theorem}[section]
\newtheorem{corollary}[theorem]{Corollary}
\newtheorem{lemma}[theorem]{Lemma}
\newtheorem{example}[theorem]{Example}

\newtheorem{remark}[theorem]{Remark}

\newcommand{\proof}{\medskip \noindent {\bf Proof. \ \ }}
\newcommand{\qed}{\null\hfill $\Box\;\;$ \medskip}

\begin{document}

\parbox{1mm}

\begin{center}
{\bf {\sc \Large On Boundedness of Calder\'on-Toeplitz Operators}}
\end{center}

\vskip 12pt

\begin{center}
{\bf Ondrej HUTN\'IK}
\end{center}


\hspace{5mm}\parbox[t]{10cm}{\fontsize{9pt}{0.1in}\selectfont\noindent{\bf
Abstract.} We study the boundedness of Toeplitz-type operators
defined in the context of the Calder\'on reproducing formula
considering the specific wavelets whose Fourier transforms are
related to Laguerre polynomials. Some sufficient conditions for
simultaneous boundedness of these Calder\'on-Toeplitz operators on
each wavelet subspace for unbounded symbols are given, where
investigating the behavior of certain sequence of iterated
integrals of symbols is helpful. A number of examples and
counterexamples is given. } \vskip 24pt

\section{Introduction}

Calder\'on-Toeplitz operators are integral operators which arise
in the context of wavelet analysis~\cite{daubechies} in connection
with the Calder\'on reproducing formula, cf.~\cite{calderon}. This
formula gives rise to a class of Hilbert spaces with reproducing
kernels, the so-called spaces of Calder\'on (or, wavelet)
transforms. These operators were formally defined
in~\cite{rochberg-CT} as a wavelet counterpart of~Toeplitz
operators defined on Hilbert spaces of~holomorphic functions.
Therefore the name "Calder\'on-Toeplitz" reflects the close
relationship with the Calder\'on reproducing formula on one side
and, on the other side, it emphasizes the fact that this operator
is unitarily equivalent to the Toeplitz-type operator
$$P_\psi M_a: W_\psi(L_2(\mathbb{R})) \to
W_\psi(L_2(\mathbb{R})),$$ where $M_a$ is the multiplication
operator by $a$ and $P_\psi$ is the orthogonal projection from
$L_2(\mathbb{R}\times \mathbb{R}_+,
v^{-2}\,\mathrm{d}u\,\mathrm{d}v)$ onto the space of wavelet
transforms
$$W_\psi(L_2(\mathbb{R}))=\Bigl\{W_\psi f(u,v)=\langle f,\psi_{u,v}\rangle;\,\, f\in
L_2(\mathbb{R})\Bigr\},$$ where $\psi\in L_2(\mathbb{R})$ is an
admissible wavelet and $\psi_{u,v}$, $(u,v)\in
\mathbb{R}\times\mathbb{R}_+$, are the shifted and scaled versions
of $\psi$, see Section~\ref{section2}. These operators are also a
useful localization tool which enables to localize a signal both
in time and frequency. For further details about localization
operators and wavelets transforms we refer to the Wong
book~\cite{wongbook}.

In~\cite{hutnik3} we have described the structure of the space of
Calder\'on transforms $W_\psi(L_2(\mathbb{R}))$ inside
$L_2(\mathbb{R}\times \mathbb{R}_+,
v^{-2}\,\mathrm{d}u\,\mathrm{d}v)$. This representation was
further used to study Calder\'on-Toeplitz operators acting on
spaces of Calder\'on transforms for general (admissible) wavelets
in~\cite{hutnik2}. Considering certain specific wavelets
in~\cite{hutnik} we were able to simplify these general results
and discover the interesting connection between the spaces of
Calder\'on transforms and poly-analytic Bergman spaces,
cf.~\cite{vasilevskibook}. This connection was further
investigated in~\cite{abreu} where also some interesting
applications were given. The specificity of this choice of
wavelets $\psi^{(k)}$, whose Fourier transforms are related to
Laguerre polynomials, has enabled us to study in~\cite{hutnikI} in
more detail the family of Calder\'on-Toeplitz operators
$T_a^{(k)}$ acting on wavelet subspaces $A^{(k)}$ (with parameter
$k=0,1,2,\dots$ being the degree of Laguerre polynomial $L_k$)
with symbols depending only on vertical variable in the upper
half-plane $\Pi$. In this specific case of wavelets the
corresponding Calder\'on-Toeplitz operators generalize classical
Toeplitz operators acting on the Bergman space in an interesting
way which differs from the case of Toeplitz operators acting on
the weighted Bergman spaces studied in~\cite{GKV}. On the other
hand, the classical Toeplitz operators and Calder\'on-Toeplitz
operators share many features in common.

In this paper we continue the detailed study of
Calder\'on-Toeplitz operators $T_a^{(k)}$ acting on wavelet
subspace $A^{(k)}$, which was initiated in~\cite{hutnikI}. For a
bounded symbol $a$ on $G$ the Calder\'on-Toeplitz operator
$T_a^{(k)}$ is clearly bounded on $A^{(k)}$. However, an
interesting and important feature of these operators on wavelet
subspaces is that they can be bounded for symbols that are
unbounded near the boundary. Therefore the aim of this paper is to
study in detail the boundedness properties of Calder\'on-Toeplitz
operators with such unbounded symbols and to give sufficient
conditions for their simultaneous boundedness on all wavelet
subspaces. The main tool in our study is the result, which we have
shown in~\cite{hutnikI}, that the Calder\'on-Toeplitz operator
$T_a^{(k)}$ acting on $A^{(k)}$ with a symbol $a=a(\Im\zeta)$,
$\zeta=(u,v)$, is unitarily equivalent to the multiplication
operator $\gamma_{a,k} I$ acting on $L_2(\mathbb{R}_+)$ with
$$\gamma_{a,k}(\xi) = \chi_+(\xi) \int_{\mathbb{R}_+} a\left(\frac{v}{2\xi}\right)
\,\ell_k^2(v)\,\mathrm{d}v,$$ where the functions
$\ell_k(x)=e^{-x/2}L_k(x)$ forms an orthonormal basis in
$L_2(\mathbb{R}_+)$ with $L_k(x)$ being the Laguerre polynomial of
order $k=0,1,\dots$ Thus, the boundedness of function
$\gamma_{a,k}$ is responsible also for the boundedness of operator
$T_a^{(k)}$. We will also show that for unbounded symbols
$a=a(v)$, $v\in\mathbb{R}_+$, the behavior of iterated means
{\setlength\arraycolsep{2pt}
\begin{eqnarray*}
C_a^{(1)}(v) & = & \int_0^v a(t)\,\mathrm{d}t, \\
C_a^{(m)}(v) & = & \int_0^v C_a^{(m-1)}(t)\,\mathrm{d}t,
\,\,\,m=2,3,\dots,
\end{eqnarray*}}rather than the behavior of
symbol $a$ itself plays a crucial role in the boundedness
properties. Contrary to the case of Toeplitz operators on weighted
Bergman spaces studied in~\cite{GKV} these means do not depend on
a weighted parameter $k$. In our case $k=0,1,\dots$ may be viewed
as a level of consideration of time-scale analysis of a signal. We
present a number of examples and construct wide families of
unbounded symbols for which the Calder\'on-Toeplitz operator is
not only bounded, but also belongs to the algebra of bounded
Calder\'on-Toeplitz operators generated by bounded symbols on
$\mathbb{R}_+$ having limits at the endpoints of $[0,+\infty]$.
This extends results for bounded symbols from~\cite{hutnikI} to
certain unbounded ones.

In the last Section~\ref{section3} we show how Calder\'on-Toeplitz
operators with unbounded symbols can appear as uniform limits of
Calder\'on-Toeplitz operators with bounded symbols.

\section{Preliminaries}\label{section2}

Here we briefly recall some necessary notations and results from
our previous works, mainly from~\cite{hutnikI}. As usual,
$\mathbb{R}$ ($\mathbb{C}$, $\mathbb{N}$) is the set of all real
(complex, natural) numbers,
$\overline{\mathbb{R}}=\mathbb{R}\cup\{-\infty,+\infty\}$ is the
two-point compactification of $\mathbb{R}$, and $\mathbb{R}_{+}$
is the positive half-line with $\chi_+$ being its characteristic
function. Let $L_2(G,\mathrm{d}\nu)$ be the space of all
square-integrable functions on $G$ with respect to measure
$\mathrm{d}\nu$, where $G=\{\zeta=(u,v);\,\, u\in \mathbb{R},
v>0\}$ is the locally compact ``$ax+b$''-group with the left
invariant Haar measure
$\mathrm{d}\nu(\zeta)=v^{-2}\,\mathrm{d}u\,\mathrm{d}v$. In what
follows we identify the group $G$ with the upper half-plane
$\Pi=\{\zeta=u+iv; \,u\in\mathbb{R}, v>0\}$ in the complex plane
$\mathbb{C}$. Then the square-integrable representation $\rho$ of
$G$ on $L_2(\mathbb{R})$ is given by
$$(\rho_{\zeta}f)(x)=f_\zeta(x)=\frac{1}{\sqrt{v}}\,f\left(\frac{x-u}{v}\right), \quad f\in
L_2(\mathbb{R}),$$ with $\zeta=(u,v)\in G$. The function $\psi \in
L_2(\mathbb{R})$ is called an \textit{admissible wavelet} if it
satisfies the so-called admissibility condition
$$\int_{\mathbb{R}_{+}} |\hat{\psi}(x\xi)|^{2}\frac{\mathrm{d}\xi}{\xi}=1$$
for almost every $x\in\mathbb{R}$, where $\hat{\psi}$ stands for
the unitary Fourier transform $\mathcal{F}: L_2(\mathbb{R})\to
L_2(\mathbb{R})$ given by
$$\mathcal{F}\{g\}(\xi)=\hat{g}(\xi)=\int_{\mathbb{R}} g(x)e^{-2\pi
ix\xi}\,\mathrm{d}x.$$

The \textit{Laguerre polynomials} $L^{(\alpha)}_n(x)$ of degree
$n\in\mathbb{Z}_+=\mathbb{N}\cup\{0\}$, and type $\alpha$ are
given by
$$L^{(\alpha)}_n(y)=\frac{y^{-\alpha} e^y}{n!}
\frac{d^n}{dy^n}\left(e^{-y}y^{n+\alpha}\right)=\sum_{k=0}^{n}{n+\alpha\choose
n-k} \frac{(-y)^{k}}{k!}, \quad y\in \mathbb{R}_+,$$
cf.~\cite[formula 8.970.1]{GR}. For $\alpha=0$ we simply write
$L_n(y)$. Recall that the system of functions
$$\ell_n(y)=e^{-y/2}L_n(y), \quad y\in \mathbb{R}_+,\,\,
n\in\mathbb{Z}_+,$$  forms an orthonormal basis in the space
$L_2(\mathbb{R}_+,\mathrm{d}y)$. For appropriate parameters
introduce the (non-negative) function
$$\Lambda_{p,m,n}^{(\alpha,\beta)}(x) = x^p e^{-x}
|L_m^{(\alpha)}(x) L_n^{(\beta)}(x)|, \quad x\in\mathbb{R}_+.$$
Then for each $\alpha,\beta\geq -\frac{1}{2}$, $\Re p>-1$,
$x\in\mathbb{R}_+$, and $m,n\in\mathbb{Z}_+$ we have
{\setlength\arraycolsep{2pt}
\begin{eqnarray}\label{upperbound1}
\Lambda_{p,m,n}^{(\alpha,\beta)}(x) & \leq & \sum_{i=0}^m
\sum_{j=0}^n \frac{(\alpha+1)_{m-i}}{(m-i)!\,i!}
\frac{(\beta+1)_{n-j}}{(n-j)!\,j!}\,x^{p+i+j} e^{-x},
\end{eqnarray}}cf.~\cite[Appendix]{hutnikI}, where $(x)_n =
\frac{\Gamma(x+n)}{\Gamma(x)}$ is the Pochhammer symbol. Thus,
{\setlength\arraycolsep{2pt}
\begin{eqnarray}\label{upperbound2}
\int_{\mathbb{R}_+}
\Lambda_{p,m,n}^{(\alpha,\beta)}(x)\,\mathrm{d}x & \leq &
\sum_{i=0}^m \sum_{j=0}^n \frac{(\alpha+1)_{m-i}}{(m-i)!\,i!}
\frac{(\beta+1)_{n-j}}{(n-j)!\,j!}\,\Gamma(p+i+j+1) \nonumber \\ &
:= & \textrm{const}_{p,m,n}^{(\alpha,\beta)}.
\end{eqnarray}}Recall
also the following closely related exact formula,
cf.~\cite[formula (16), p. 330]{WG}, {\setlength\arraycolsep{2pt}
\begin{eqnarray}\label{Laguerreformula}
& \phantom{=} & \int_{\mathbb{R}_+} x^{p} e^{-x}
L_m^{(\alpha)}(x)L_n^{(\beta)}(x)\,\mathrm{d}x \nonumber\\ & = &
\Gamma(p+1) \sum_{i=0}^{\min\{m,n\}} (-1)^{m+n}{p-\alpha \choose
m-i} {p-\beta \choose n-i} {p+i \choose i},
\end{eqnarray}}where $\Re p>-1$, $\alpha,\beta>-1$,
$m,n\in\mathbb{Z}_+$, and $${a \choose b} =
\frac{\Gamma(a+1)}{\Gamma(b+1) \Gamma(a-b+1)}.$$

Further, for $k\in\mathbb{Z}_+$ we consider the functions
(admissible wavelets) $\psi^{(k)}$ and $\bar{\psi}^{(k)}$ on
$\mathbb{R}$, whose Fourier transforms are given by
$$\hat{\psi}^{(k)}(\xi)=\chi_+(\xi)\sqrt{2\xi}\,\ell_k(2\xi)\quad \textrm{and}
\quad \hat{\bar{\psi}}^{(k)}(\xi)=\hat{\psi}^{(k)}(-\xi),$$
respectively. Let $A^{(k)}$, resp. $\bar{A}^{(k)}$, be the spaces
of wavelet transforms of functions $f\in H^+_2(\mathbb{R})$, resp.
$f\in H^-_2(\mathbb{R})$, with respect to wavelets $\psi^{(k)}$,
resp. $\bar{\psi}^{(k)}$, where {\setlength\arraycolsep{2pt}
\begin{eqnarray*}
H^{+}_2(\mathbb{R}) & = & \left\{f\in L_2(\mathbb{R}); \,\,\textrm{supp}\, \hat{f}\subseteq [0,+\infty)\right\}, \\
H^-_2(\mathbb{R}) & = & \left\{f\in L_2(\mathbb{R});
\,\,\textrm{supp}\, \hat{f}\subseteq (-\infty,0]\right\},
\end{eqnarray*}}are the Hardy spaces, respectively.

Consider the unitary operators $$U_{1}= (\mathcal{F}\otimes I):
L_2(G, \mathrm{d}\nu(\zeta)) \to L_{2}(\mathbb{R}, \mathrm{d}u)
\otimes L_{2}(\mathbb{R}_{+}, v^{-2}\mathrm{d}v)$$ with $\zeta =
(u,v) \in G$, and
$$U_{2}: L_{2}(\mathbb{R},\,\mathrm{d}u)\otimes
L_{2}(\mathbb{R}_{+},\,v^{-2}\mathrm{d}v) \to
L_{2}(\mathbb{R},\,\mathrm{d}x)\otimes
L_{2}(\mathbb{R}_{+},\,\mathrm{d}y)$$ given by
$$U_{2}: F(u,v) \mapsto
\frac{\sqrt{2|x|}}{y} F\left(x,\frac{y}{2|x|}\right).$$
In~\cite{hutnik} we have proved the following important result
describing the structure of wavelet subspaces $A^{(k)}$ inside
$L_2(G,\mathrm{d}\nu)$.

\begin{theorem}\label{thm1}
The unitary operator $U=U_{2}U_{1}$ gives an isometrical
isomorphism of the space $L_{2}(G,\mathrm{d}\nu)$ onto
$L_{2}(\mathbb{R},\,\mathrm{d}x)\otimes
L_{2}(\mathbb{R}_{+},\,\mathrm{d}y)$ under which the wavelet
subspace $A^{(k)}$ is mapped onto $L_2(\mathbb{R}_+) \otimes
L_{k}$, where $L_{k}$ is the rank-one space generated by function
$\ell_{k}(y)=e^{-y/2}L_k(y)$.
\end{theorem}

This result is a "wavelet" analog of results obtained for the
Bergman and poly-Bergman spaces, cf.~\cite{vasilevski2}, and
enables to study an interesting connection between wavelet spaces
related to Laguerre polynomials and poly-Bergman spaces in more
detail, which is done in paper~\cite{abreu}.

Following the general scheme presented in~\cite{vasilevski}, let
us introduce the isometric imbedding
$$Q_{k}: L_{2}(\mathbb{R}_+) \to
L_{2}(\mathbb{R})\otimes L_{2}(\mathbb{R}_{+})$$ given by
$$\left(Q_k f\right)(x,y) = \chi_+(x) f(x)\ell_{k}(y),$$ here the function $f$
is extended to an element of $L_2(\mathbb{R})$ by setting
$f(x)\equiv 0$ for $x<0$. The adjoint operator
$$Q_k^{*}: L_{2}(\mathbb{R})\otimes
L_{2}(\mathbb{R}_{+}) \to L_{2}(\mathbb{R}_+)$$ is given by
$$\left(Q_k^{*} F\right)(x) = \chi_{+}(x) \int_{\mathbb{R}_{+}}
F(x,\tau)\ell_{k}(\tau)\,\mathrm{d}\tau,$$ and we have

\begin{theorem}\label{IIR*}
The operator $R_{k}: A^{(k)} \to L_{2}(\mathbb{R}_+)$, where
\begin{equation}\label{RF}
(R_k F)(\xi) = \chi_+(\xi)
\sqrt{2\xi}\int_{\mathbb{R}\times\mathbb{R}_+} F(u,v)
\ell_k(2v\xi) e^{-2\pi i\xi u}\,\frac{\mathrm{d}u\mathrm{d}v}{v},
\end{equation} is an
isometrical isomorphism admitting decomposition $R_k=Q_k^*
U_2(\mathcal{F}\otimes I)$.
\end{theorem}

\begin{corollary}\rm
The inverse isomorphism $R_{k}^*: L_2(\mathbb{R}_+)\to A^{(k)}$
given by
\begin{equation}\label{R*}
\left(R_{k}^{*}f\right)(u,v) = \sqrt{2}v \int_{\mathbb{R}_+}
f(\xi)\ell_k(2\xi v)e^{2\pi i\xi
u}\,\sqrt{\xi}\,\mathrm{d}\xi\end{equation} admits the
decomposition $R_k^* = (\mathcal{F}^{-1}\otimes I) U_2^{-1} Q_k$
with $\mathcal{F}^{-1}: L_2(\mathbb{R})\to L_2(\mathbb{R})$ being
the inverse Fourier transform.
\end{corollary}

The above representation of wavelet subspaces is especially
important in the study of Toeplitz-type operators related to
wavelets whose symbols depend only on vertical variable $v=\Im
\zeta$ in the upper half-plane $\Pi$ of the complex plane
$\mathbb{C}$. For a given $L_\infty(G,\mathrm{d}\nu)$-function
$a(\zeta)=a(v)$ depending only on $v = \Im\zeta$, $\zeta\in G$,
define the \textit{Calder\'on-Toeplitz operator} $T_a^{(k)}:
A^{(k)}\to A^{(k)}$ with symbol $a$ as
$$T_{a}^{(k)} = P^{(k)}M_{a},$$ where $M_a$ is the operator
of pointwise multiplication by $a$ and $P^{(k)}$ is the orthogonal
projection from $L_2(G,\mathrm{d}\nu)$ to the wavelet subspace
$A^{(k)}$. Formally, these operators were introduced
in~\cite{rochberg-CT} and further studied in general e.g. in
papers~\cite{hutnik2}, \cite{nowak3}, \cite{nowak4}
and~\cite{rochberg3}. The following important result, which
enables to reduce Calder\'on-Toeplitz operator to a certain
multiplication operator, was proved in~\cite{hutnikI}. In fact, it
is the main tool of our study.

\begin{theorem}\label{CTO1}
Let $a=a(v)$, $v\in\mathbb{R}_+$, be a measurable symbol on $G$.
Then the Calder\'on-Toeplitz operator $T_{a}^{(k)}$ acting on
$A^{(k)}$ is unitarily equivalent to the multiplication operator
$$\gamma_{a,k} I = R_{k} T_{a}^{(k)} R^{*}_k$$ acting on
$L_2(\mathbb{R}_+)$, where $R_k$ and $R^*_k$ are given
by~$\mathrm{(\ref{RF})}$ and~$\mathrm{(\ref{R*})}$, respectively.
The function $\gamma_{a,k}$ is given by
\begin{equation}\label{gamma1}
\gamma_{a,k}(\xi) =
\int_{\mathbb{R}_{+}}a\left(\frac{v}{2\xi}\right)\ell_k^2(v)\,\mathrm{d}v,
\quad \xi\in\mathbb{R}_+.
\end{equation}
\end{theorem}

Note that all the above results and definitions may be stated
analogously for the space $\bar{A}^{(k)}$. In what follows we
restrict our attention only to wavelet subspaces $A^{(k)}$ and
operators $T_a^{(k)}$ acting on them. For more results on
properties of Calder\'on-Toeplitz operators $T_a^{(k)}$ with
symbols $a=a(v)$ and properties of the corresponding function
$\gamma_{a,k}(\xi)$ responsible for many interesting features of
these operators and their algebras, see the recent
paper~\cite{hutnikI}.

\section{Boundedness of Calder\'on-Toeplitz
operator}\label{section5}

Clearly, if $a=a(v)$ is a bounded symbol on $G$, then the operator
$T_a^{(k)}$ is bounded on $A^{(k)}$, and for its operator norm
holds
$$\|T_a^{(k)}\| \leq \textrm{ess-sup}\,|a(v)|.$$ Thus all spaces
$A^{(k)}$, $k\in\mathbb{Z}_+$, are naturally appropriate for
Calder\'on-Toeplitz operators with bounded symbols. However, we
may observe that the result of Theorem~\ref{CTO1} suggests
considering not only $L_{\infty}(G,\mathrm{d}\nu)$-symbols, but
also \textit{unbounded} ones. In this case we obviously have

\begin{corollary}\rm\label{corgamma}
Calder\'on-Toeplitz operator $T_a^{(k)}$ with a measurable symbol
$a=a(v)$, $v\in\mathbb{R}_+$, is bounded on $A^{(k)}$ if and only
if the function $\gamma_{a,k}(\xi)$ is bounded on $\mathbb{R}_+$,
and
$$\|T_a^{(k)}\| = \sup_{\xi \in\mathbb{R}_+}
|\gamma_{a,k}(\xi)|.$$
\end{corollary}

From this result we immediately have that the Calder\'on-Toeplitz
operator $T_a^{(1)}$ with unbounded symbol
\begin{equation}\label{exmunbounded3}
a(v)=\frac{1}{\sqrt{v}}\sin\frac{1}{v}, \quad v\in\mathbb{R}_+,
\end{equation} \textit{is bounded} on $A^{(1)}$
because the corresponding function
$$\gamma_{a,1}(\xi) = \frac{\sqrt{2\pi}}{4} e^{-2\sqrt{\xi}}
\left[(2\sqrt{\xi}-8\xi)\frac{\cos
2\sqrt{\xi}}{2\sqrt{\xi}}+(3-2\sqrt{\xi})\frac{\sin
2\sqrt{\xi}}{2\sqrt{\xi}}\right], \quad
\xi\in\overline{\mathbb{R}}_+,$$ is bounded, see~\cite[Example
4.4]{hutnikI}. However, due to computational limitations (to find
an explicit form of the function $\gamma_{a,k}(\xi)$) we can not
say anything about the boundedness of $T_a^{(k)}$ for arbitrary
$k$. Fortunately, according to Theorem~\ref{thmmeans1} we will be
able to show much more for a more general class of unbounded
symbols including that symbol given by~(\ref{exmunbounded3}).

\begin{example}\rm\label{exmunbounded}
For oscillating symbol $a(v)=e^{2v \mathrm{i}}$ (with
$\mathrm{i}^2=-1$) we have
$$\gamma_{a,k}(\xi) =
\frac{(-1)^k}{(\xi-\mathrm{i})^{2k+1}}\sum_{j=0}^{k}
(-1)^{j}\left[{k\choose j}\right]^2 \xi^{2j+1}, \quad
\xi\in\overline{\mathbb{R}}_+,$$ see~\cite[Example 4.5]{hutnikI}.
Thus, the Calder\'on-Toeplitz operator $T_a^{(k)}$ acting on
$A^{(k)}$ is bounded for each $k\in\mathbb{Z}_+$, and moreover
$\gamma_{a,k}(\xi)\in C[0,+\infty]$.
\end{example}

In both the above mentioned examples (more precisely, in the first
one just for the case $k=1$, but we will show later that also for
all $k\in\mathbb{Z}_+$) we are in situation that the
Calder\'on-Toeplitz operator $T_a^{(k)}$ belongs to the
$C^*$-algebra
$\mathcal{T}_k\left(L_\infty^{\{0,+\infty\}}(\mathbb{R}_+)\right)$
generated by (bounded) Calder\'on-Toeplitz operators with
$L_\infty(\mathbb{R}_+)$-\-sym\-bols having limits at the points
$0$ and $+\infty$, cf.~\cite[Section 4]{hutnikI}. Recall that the
Calder\'on-Toeplitz operator $T_a^{(k)}$ with a symbol $a=a(v)$
belongs to the algebra
$\mathcal{T}_k\left(L_\infty^{\{0,+\infty\}}(\mathbb{R}_+)\right)$
if and only if the corresponding function $\gamma_{a,k}(\xi)$
belongs to $C[0,+\infty]$. This means that the algebra
$\mathcal{T}_k\left(L_\infty^{\{0,+\infty\}}(\mathbb{R}_+)\right)$
contains many more Calder\'on-Toeplitz operators than was
described in~\cite{hutnikI}, because it also contains (bounded)
Calder\'on-Toeplitz operators whose (generally unbounded) symbols
$a(v)$ need not have limits at the endpoints $0$ and $+\infty$.

\begin{example}\rm
An easy example of unbounded symbol, for which the
Calde\-r\'on-Toeplitz operator is unbounded for each
$k\in\mathbb{Z}_+$, is the function $a=a(v)=v^p$ with $\Re p>-1$
and $\Re p\neq 0$. The explicit formula for $\gamma_{a,k}$ has the
form
$$\gamma_{a,k}(\xi) = \int_{\mathbb{R}_+} a\left(\frac{v}{2\xi}\right)
\ell_k^2(v)\,\mathrm{d}v = \frac{1}{(2\xi)^p} \int_{\mathbb{R}_+}
\Lambda_{p,k,k}^{(0,0)}(v)\,\mathrm{d}v, \quad \xi\in\mathbb{R}_+.
$$ Since by the formula~(\ref{Laguerreformula}) the last integral is a
constant, the function $\gamma_{a,k}(\xi)$ is clearly unbounded on
$\mathbb{R}_+$. Thus the operator $T_a^{(k)}$ is not bounded on
$A^{(k)}$ and does not belong to the algebra
$\mathcal{T}_k\left(L_\infty^{\{0,+\infty\}}(\mathbb{R}_+)\right)$.
This case of symbols is subsumed in Theorem~\ref{thm3}.
\end{example}


As we have shown in~\cite[Theorem 4.2]{hutnikI}, the behavior of a
bounded function $a(v)$ near the point $0$, or $+\infty$
determines the behavior of function $\gamma_{a,k}(\xi)$ near the
point $+\infty$, or $0$, respectively. The existence of limits of
$a(v)$ in these endpoints guarantees the continuity of
$\gamma_{a,k}(\xi)$ on the whole $\overline{\mathbb{R}}_+$,
however this condition is not necessary even for bounded symbols,
see~\cite[Remark 4.3]{hutnikI}. Continuity of function
$\gamma_{a,k}$ on the whole $\overline{\mathbb{R}}_+$ then
guarantees its boundedness, and therefore by
Corollary~\ref{corgamma} the boundedness of the corresponding
Calder\'on-Toeplitz operator $T_a^{(k)}$ on wavelet subspace
$A^{(k)}$.

However, as example of symbol~(\ref{exmunbounded3}) shows (we will
do it exactly and more generally in Example~\ref{exmunbounded2})
the Calder\'on-Toeplitz operator $T_a^{(k)}$ can be
\textit{bounded} and \textit{belong to the algebra
$\mathcal{T}_k\left(L_\infty^{\{0,+\infty\}}(\mathbb{R}_+)\right)$
for each} $k\in\mathbb{Z}_+$ even for \textit{unbounded symbols}
$a=a(v)$. Thus, in what follows we study this phenomena in more
detail and will be interested in unbounded symbols to have a
sufficiently large class of them common to all admissible $k$. For
this purpose denote by $L_{1}(\mathbb{R}_+, 0)$ the class of
functions $a=a(v)$ such that
$$a(v)e^{-\varepsilon v} \in L_1(\mathbb{R}_+), \,\,\,\textrm{for
any}\,\,\, \varepsilon > 0.$$ We give some conditions on the
behavior of $L_{1}(\mathbb{R}_+, 0)$-symbols (in fact on the
behavior of certain means of these symbols) which guarantees the
boundedness of function $\gamma_{a,k}(\xi)$.

\begin{remark}\rm
Using several formulas, more precisely~\cite[formula 8.976.3]{GR},
\cite[formula (5), p. 209]{Rainville} and~\cite[formula
8.976.1]{GR}, we may rewrite the function $\gamma_{a,k}(\xi)$ as
follows \small{{\setlength\arraycolsep{2pt}
\begin{eqnarray*}
& \phantom{=} & \gamma_{a,k}(\xi) \\ & = & \frac{1}{2^{2k+1}}
\sum_{i=0}^k \sum_{j=0}^{2k} \sum_{r=0}^{j} {2k-2i \choose k-i}
{2i \choose j} {j\choose r} {2i\choose i}\frac{(-1)^r }{r!}
(1-4\xi)^{2i-j} (4\xi)^{j+1} I_r(\xi),
\end{eqnarray*}}}\normalsize where $$I_r(\xi) = \int_{\mathbb{R}_+} a(v)v^r
e^{-2v\xi}\,\mathrm{d}v.$$ Then the last integral is, in fact, the
integral in the formula of function $\gamma_{a,\lambda}(x)$ for
Toeplitz operators on the upper half-plane (the so-called
parabolic case), see~\cite[formula (13.1.1), p.
329]{vasilevskibook}, or~\cite[formula (2.6)]{GKV}. Therefore it
would be natural to consider certain means of symbols depending on
parameter $k$ as it is done therein, but we will not do it this
way and our means of symbols will not depend on weight parameter.
\end{remark}

For any $L_{1}(\mathbb{R}_+, 0)$-symbol $a(v)$ define the
following averaging functions {\setlength\arraycolsep{2pt}
\begin{eqnarray*}
C_a^{(1)}(v) & = & \int_0^v a(t)\,\mathrm{d}t, \\
C_a^{(m)}(v) & = & \int_0^v C_a^{(m-1)}(t)\,\mathrm{d}t,
\,\,\,m=2,3,\dots
\end{eqnarray*}}The functions $C_a^{(m)}$ constitute a "sequence of iterated
integrals" of symbol $a$.

\begin{theorem}\label{thmmeans1}
Let $a=a(v)\in L_{1}(\mathbb{R}_+, 0)$ and for any
$m\in\mathbb{N}$ suppose that the function $C_a^{(m)}$ has the
following asymptotic behavior
\begin{equation}\label{asymptoticC1}
C_a^{(m)}(v) = \mathcal{O}(v^m),\,\,\,\, \textrm{as}\,\,\,\, v\to
0,
\end{equation} and
\begin{equation}\label{asymptoticC2}
C_a^{(m)}(v) = \mathcal{O}(v^m),\,\,\,\, \textrm{as}\,\,\,\, v\to
+\infty.
\end{equation} Then for each $k\in\mathbb{Z}_+$ we have
$$\sup_{\xi\in\mathbb{R}_+} |\gamma_{a,k}(\xi)|<+\infty.$$ Consequently, the
corresponding Calder\'on-Toeplitz operator $T_a^{(k)}$ is bounded
on $A^{(k)}$ for every $k\in\mathbb{Z}_+$.
\end{theorem}

\proof Let $m\geq 1$. The condition~(\ref{asymptoticC1}) together
with the condition~(\ref{asymptoticC2}) imply that for all
$v\in\mathbb{R}_+$ the estimate
\begin{equation}\label{asymptoticC3}
|C_a^{(m)}(v)|\leq \textrm{const}\,v^m
\end{equation} holds, where "$\textrm{const}$" does not depend on $v\in\mathbb{R}_+$.
Integrating by parts $m$-times we have for all
$\xi\in\mathbb{R}_+$ that {\setlength\arraycolsep{2pt}
\begin{eqnarray*}
\gamma_{a,k}(\xi) & = & 2\xi \int_{\mathbb{R}_+} \ell_k^2(2v\xi)
\,\mathrm{d}C_a^{(1)}(v) \\ & = & -2\xi \int_{\mathbb{R}_+}
\frac{\mathrm{d}}{\mathrm{d}v}
\ell_k^2(2v\xi)\,\mathrm{d}C_a^{(2)}(v) \\ & \phantom{=} & \hspace{2cm} \vdots \\
& = & (-1)^m 2\xi \int_{\mathbb{R}_+}
\frac{\mathrm{d}^m}{\mathrm{d}v^m}\ell_k^2(2v\xi)\,\mathrm{d}C_a^{(m+1)}(v)
\\ & = & (-1)^m 2\xi \int_{\mathbb{R}_+} C_a^{(m)}(v)
\frac{\mathrm{d}^m}{\mathrm{d}v^m}\ell_k^2(2v\xi)\,\mathrm{d}v.
\end{eqnarray*}}Using the estimates~(\ref{asymptoticC3})
and~(\ref{upperbound2}) we then get {\setlength\arraycolsep{2pt}
\begin{eqnarray*}
& \phantom{=} & |\gamma_{a,k}(\xi)| \\ & \leq & (2\xi)^{m+1}
\sum_{i=0}^m \sum_{j=0}^i {m \choose i} {i \choose j}
\int_{\mathbb{R}_+} |C_a^{(m)}(v)| e^{-2v\xi}
|L_{k-i+j}^{(i-j)}(2v\xi) L_{k-j}^{(j)}(2v\xi)|\,\mathrm{d}v \\ &
\leq & \textrm{const}\, \sum_{i=0}^m \sum_{j=0}^i {m \choose i} {i
\choose j} \int_{\mathbb{R}_+} (2\xi v)^m e^{-2v\xi}
|L_{k-i+j}^{(i-j)}(2v\xi) L_{k-j}^{(j)}(2v\xi)|\,2\xi\,\mathrm{d}v \\
& = & \textrm{const}\,\sum_{i=0}^m \sum_{j=0}^i {m \choose i} {i
\choose j} \int_{\mathbb{R}_+}
\Lambda_{m,k-i+j,k-j}^{(i-j,j)}(x)\,\mathrm{d}x \\ & \leq &
\textrm{const}\,\sum_{i=0}^m \sum_{j=0}^i {m \choose i} {i \choose
j} \,\textrm{const}_{m,k-i+j,k-j}^{(i-j,j)} \\ & < & +\infty,
\end{eqnarray*}}thus the function $\gamma_{a,k}$ is bounded for each
$k\in\mathbb{Z}_+$, and by Corollary~\ref{corgamma} the
Calder\'on-Toeplitz operator $T_a^{(k)}$ is bounded as well. \qed

\begin{remark}\rm
The condition~(\ref{asymptoticC1}) guarantees the boundedness of
the function $\gamma_{a,k}(\xi)$ at a neighborhood of
$\xi=+\infty$, while the condition~(\ref{asymptoticC2}) guarantees
its boundedness at a neighborhood of $\xi=0$. Observe that if the
conditions~(\ref{asymptoticC1}) and~(\ref{asymptoticC2}) hold for
some $m=m_0$, then according to~(\ref{asymptoticC3}) hold also for
$m=m_0+1$. Indeed, $$|C_a^{(m_0+1)}(v)| \leq \int_0^v
|C_a^{(m_0)}(t)|\,\mathrm{d}t \leq \textrm{const}\,\int_0^v
t^{m_0}\,\mathrm{d}t \leq \textrm{const}\,v^{m_0+1}.$$
\end{remark}

The main advantage of Theorem~\ref{thmmeans1} is that we need not
have an explicit form of the corresponding function $\gamma_{a,k}$
for an unbounded symbol $a=a(v)$ to decide about its boundedness,
as it is e.g. in Example~\ref{exmunbounded}.

\begin{theorem}\label{thmmeans2}
Let $a=a(v)\in L_1(\mathbb{R}_+, 0)$. If for any
$m,n\in\mathbb{N}$, any $\lambda_1\in\mathbb{R}_+$ and any
$\lambda_2\in (0,n+1)$ holds
\begin{equation}\label{asymptoticC4}
C_a^{(m)}(v) = \mathcal{O}\left(v^{m+\lambda_1}\right),\,\,\,\,
\textrm{as}\,\,\,\, v\to 0, \end{equation} and
\begin{equation}\label{asymptoticC5}
C_a^{(n)}(v) = \mathcal{O}\left(v^{n-\lambda_2}\right),\,\,\,\,
\textrm{as}\,\,\,\, v\to +\infty,
\end{equation} then for each $k\in\mathbb{Z}_+$ we have
$$\lim_{\xi\to+\infty} \gamma_{a,k}(\xi) = 0 = \lim_{\xi\to 0} \gamma_{a,k}(\xi).$$
\end{theorem}

\proof Analogously as in the proof of Theorem~\ref{thmmeans1}
using the condition~(\ref{asymptoticC4}) we have
$$|\gamma_{a,k}(\xi)| \leq \frac{1}{(2\xi)^{\lambda_1}}\,
\textrm{const}\,\sum_{i=0}^m \sum_{j=0}^i {m \choose i} {i \choose
j} \int_{\mathbb{R}_+}
\Lambda_{m+\lambda_1,k-i+j,k-j}^{(i-j,j)}(x)\,\mathrm{d}x,$$ where
"$\textrm{const}$" does not depend on $v\in\mathbb{R}_+$. Thus,
$\lim\limits_{\xi\to+\infty} \gamma_{a,k}(\xi) = 0$.

Similarly, integrating by parts $n$-times and using the
condition~(\ref{asymptoticC5}) we get
$$|\gamma_{a,k}(\xi)| \leq (2\xi)^{\lambda_2}\,
\textrm{const}\,\sum_{i=0}^n \sum_{j=0}^i {n \choose i} {i \choose
j} \int_{\mathbb{R}_+}
\Lambda_{n-\lambda_2,k-i+j,k-j}^{(i-j,j)}(x)\,\mathrm{d}x,$$ where
(the different) "$\textrm{const}$" does not depend on
$v\in\mathbb{R}_+$. Letting $\xi\to 0$ we have again the desired
result. \qed

In other words, Theorem~\ref{thmmeans2} gives the condition on the
behavior of $L_1(\mathbb{R}_+,0)$-\-sym\-bols such that the
function $\gamma_{a,k}(\xi)\in C[0,+\infty]$, and thus the
corresponding Calder\'on-Toeplitz operator $T_a^{(k)}$ belongs to
the algebra
$\mathcal{T}_k\left(L_\infty^{\{0,+\infty\}}(\mathbb{R}_+)\right)$.
In the next example we present a wide class of oscillating symbols
$a=a(v)$ for which the Calder\'on-Toeplitz operator $T_a^{(k)}$
belongs to
$\mathcal{T}_k\left(L_\infty^{\{0,+\infty\}}(\mathbb{R}_+)\right)$
for each $k\in\mathbb{Z}_+$.

\begin{example}\rm\label{exmunbounded2}
For $\alpha>0$ and $\beta \in (0,1)$ consider the unbounded symbol
$$a(v)=v^{-\beta}\sin v^{-\alpha}, \quad v\in\mathbb{R}_+.$$ However, the function $a(v)$ is
continuous at $v=+\infty$ for all admissible values of parameters,
and therefore $\gamma_{a,k}(0) = a(+\infty) = 0$. On the other
side, it is difficult to verify the behavior of function
$\gamma_{a,k}(\xi)$ at the endpoint $+\infty$ by a direct
computation. According to~\cite[Example~13.1.4]{vasilevskibook} we
have
$$C_a^{(1)}(v) = \frac{v^{\alpha-\beta+1}}{\alpha} \cos v^{-\alpha}
+ \mathcal{O}(v^{2\alpha-\beta+1}), \quad \textrm{as}\,\,\, v\to
0.$$
From it follows that for $\alpha>\beta$ the first condition
in~(\ref{asymptoticC4}) holds for $m=1$ and
$\lambda_1=\alpha-\beta$. By Theorem~\ref{thmmeans2} the function
$\gamma_{a,k}(\xi)$ is bounded, and therefore the corresponding
Calder\'on-Toeplitz operator $T_a^{(k)}$ is bounded for each $k\in
\mathbb{Z}_+$. Observe that this is in accordance with the
obtained result of a special case for the
symbol~(\ref{exmunbounded3}) and $k=1$. Here we have extended it
for much more general class of unbounded symbols and the whole
range of parameters $k$.

If $\alpha\leq\beta$, then
$$C_a^{(m)}(v) =
\mathcal{O}(v^{m\alpha-\beta+m}),\,\,\,\,\textrm{as}\,\,\,\, v\to
0.$$ Thus for each $\alpha\leq\beta$ there exists
$m_0\in\mathbb{N}$ such that $m_0\alpha>\beta$, and therefore the
first condition in~(\ref{asymptoticC4}) holds for $m=m_0$ and
$\lambda_1=m_0\alpha-\beta$, which guarantees that
$\gamma_{a,k}(\xi)$ is continuous at $\xi=0$. Since for all
parameters $\alpha>0$ and $\beta\in (0,1)$ the function
$\gamma_{a,k}(\xi)\in C[0,+\infty]$, then each Calder\'on-Toeplitz
operator $T_a^{(k)}$ belongs to the algebra
$\mathcal{T}_k\left(L_\infty^{\{0,+\infty\}}(\mathbb{R}_+)\right)$
for each $k\in\mathbb{Z}_+$.
\end{example}

In fact, Theorem~\ref{thmmeans2} partially extends the
result~\cite[Theorem 4.2]{hutnikI} stated for bounded symbols to
certain unbounded ones. In Example~\ref{exmunbounded} and
Example~\ref{exmunbounded2} we have provided such oscillating
symbols $a=a(v)$ for which the Calder\'on-Toeplitz operator
$T_a^{(k)}$ belongs to the algebra
$\mathcal{T}_k\left(L_\infty^{\{0,+\infty\}}(\mathbb{R}_+)\right)$
for the whole range of parameters $k$. Now we give an example of a
bounded oscillating symbol such that the bounded operator
$T_a^{(k)}$ does not belong to the algebra
$\mathcal{T}_k\left(L_\infty^{\{0,+\infty\}}(\mathbb{R}_+)\right)$.

\begin{example}\rm The function
$$a(v) = v^\mathrm{i} = e^{\mathrm{i} \ln v}, \quad v\in\mathbb{R}_+,
$$ is oscillating near the endpoints $0$ and $+\infty$, but it is
bounded on $\mathbb{R}_+$, and therefore the Calder\'on-Toeplitz
operator $T_a^{(k)}$ is bounded for each $k\in\mathbb{R}_+$.
Changing the variable $x=2v\xi$ yields
{\setlength\arraycolsep{2pt}
\begin{eqnarray*}
\gamma_{a,k}(\xi) & = & 2\xi\int_{\mathbb{R}_+} v^\mathrm{i}
\ell_k^2(2v\xi)\,\mathrm{d}v = (2\xi)^{-\mathrm{i}}
\int_{\mathbb{R}_+} x^\mathrm{i} \ell_k^2(x)\,\mathrm{d}x \\ & = &
(2\xi)^{-\mathrm{i}} \int_{\mathbb{R}_+}
\Lambda_{\mathrm{i},k,k}^{(0,0)}(x)\,\mathrm{d}x.
\end{eqnarray*}}Since by the formula~(\ref{Laguerreformula}) the last
integral is a constant depending on $k$, the function
$\gamma_{a,k}(\xi)$ oscillates and has no limit when $\xi\to 0$ as
well as when $\xi\to +\infty$. Thus the bounded
Calder\'on-Toeplitz operator $T_a^{(k)}$ \textit{does not belong}
to the algebra
$\mathcal{T}_k\left(L_\infty^{\{0,+\infty\}}(\mathbb{R}_+)\right)$.
Hence not all oscillating symbols (even bounded and continuous)
generate an operator from
$\mathcal{T}_k\left(L_\infty^{\{0,+\infty\}}(\mathbb{R}_+)\right)$.
\end{example}


In the following theorem we show that \textit{the boundedness of a
Toeplitz operator on the Bergman space} with non-negativity of
symbol or its means guarantees the boundedness of
Calder\'on-Toeplitz operator on each wavelet subspace as well.

\begin{theorem}
(i) Let $a=a(v)\in L_{1}(\mathbb{R}_+, 0)$ be non-negative almost
everywhere. If $T_a^{(0)}$ is bounded on $A^{(0)}$, then the
operator $T_a^{(k)}$ is bounded on $A^{(k)}$ for each
$k\in\mathbb{Z}_+$.

(ii) Let $C_a^{(m)}$ be non-negative almost everywhere for a
certain $m=m_0$. If $T_a^{(0)}$ is bounded on $A^{(0)}$, then the
operator $T_a^{(k)}$ is bounded on $A^{(k)}$ for each
$k\in\mathbb{Z}_+$.
\end{theorem}

\proof (i) From assumptions we have $$\gamma_{a,0}(\xi) = 2\xi
\int_{\mathbb{R}_+} a(v) e^{-2v\xi}\,\mathrm{d}v \geq 2\xi
\int_{0}^{(2\xi)^{-1}} a(v) e^{-2v\xi}\,\mathrm{d}v \geq
\frac{2\xi}{e}\,C_a^{(1)}\left((2\xi)^{-1}\right).$$ Putting
$(2\xi)^{-1} = v$ we get
$$C_a^{(1)}(v) \leq \left(e \sup_{\xi\in\mathbb{R}_+} |\gamma_{a,0}(\xi)|\right)\, v = \textrm{const}\,v,$$
which by Theorem~\ref{thmmeans1} means that $T_a^{(k)}$ is bounded
on $A^{(k)}$ for each $k\in\mathbb{Z}_+$.

(ii) Integrating by parts $m_0$-times we obtain
{\setlength\arraycolsep{2pt}
\begin{eqnarray*}
\gamma_{a,0}(\xi) & = & (-1)^{m_0}\,2\xi \int_{\mathbb{R}_+}
C_a^{(m_0)}(v) \frac{\mathrm{d}^{m_0}}{\mathrm{d}v^{m_0}} e^{-2v\xi}\,\mathrm{d}v \\
& = &
(2\xi)^{m_0+1}\int_{\mathbb{R}_+} C_a^{(m_0)}(v) e^{-2v\xi}\,\mathrm{d}v \\
& \geq & (2\xi)^{m_0+1} \int_{0}^{(2\xi)^{-1}} C_a^{(m_0)}(v)
e^{-2v\xi}\,\mathrm{d}v \\ & \geq & (2\xi)^{m_0+1} e^{-1}
C_a^{(m_0+1)}\left((2\xi)^{-1}\right).
\end{eqnarray*}}Again putting $(2\xi)^{-1} = v$ we have
$$C_a^{(m_0+1)}(v) \leq \left(e \sup_{\xi\in\mathbb{R}_+} |\gamma_{a,0}(\xi)|\right)\, v^{m_0+1}$$
and by Theorem~\ref{thmmeans1} the boundedness of $T_a^{(k)}$ on
$A^{(k)}$ for each $k\in\mathbb{Z}_+$ follows. \qed

According to the presented examples an unbounded symbol must have
a sufficiently sophisticated oscillating behavior at neighborhoods
of the points $0$ and $+\infty$ to generate a bounded
Calder\'on-Toeplitz operator. In what follows we show that
infinitely growing positive symbols (as in the case of identity,
or its powers) cannot generate bounded Calder\'on-Toeplitz
operators in general. For this purpose for a non-negative function
$a=a(v)$ put
$$\theta_{a}(v)=\inf_{t\in (0,v)} a(t)\,\,\,\,\textrm{and}\,\,\,\,\Theta_{a}(v)=\inf_{t\in (v/2,v)} a(t).$$

\begin{theorem}\label{thm3}
For a given non-negative symbol $a=a(v)$ if either
\begin{equation}\label{limm_a,0}
\lim_{v\to 0} \theta_a(v)=+\infty
\end{equation} or
\begin{equation}\label{limm_a,infty}
\lim_{v\to +\infty} \Theta_{a}(v)=+\infty,
\end{equation} then the Calder\'on-Toeplitz operator $T_a^{(k)}$
is unbounded on each $A^{(k)}$, $k\in\mathbb{Z}_+$.
\end{theorem}

\proof If the condition~(\ref{limm_a,0}) holds, then
$$C_a^{(1)}(v) = \int_0^v a(t)\,\mathrm{d}t \geq v\, \theta_a(v),$$ which yields
$v^{-1} C_a^{(1)}(v) \to +\infty$, as $v\to 0$.

If the condition~(\ref{limm_a,infty}) holds, then
$$v^{-1}C_a^{(1)}(v) > v^{-1} \int_{v/2}^v a(t)\,\mathrm{d}t \geq \frac{1}{2} \Theta_{a}(v),$$
which again yields $v^{-1} C_a^{(1)}(v) \to +\infty$, as $v\to
+\infty$. \qed

\begin{example}\rm
For the family of non-negative symbols on $\mathbb{R}_+$ in the
form
$$a(v)=v^{-\beta}\ln^2 v^{-\alpha}, \quad \beta\in[0,1],\,
\alpha>0,$$ we have that for all admissible parameters holds
$\lim\limits_{v\to 0} \theta_a(v) = +\infty$, and thus by
Theorem~\ref{thm3} the Calder\'on-Toeplitz operator $T_a^{(k)}$ is
unbounded on $A^{(k)}$ for each $k\in\mathbb{Z}_+$.
\end{example}

In the following example we use the result of Theorem~\ref{thm3}
to study boundedness of a Calder\'on-Toeplitz operator with
unbounded symbol as a product of two symbols for which the
corresponding Calder\'on-Toeplitz operators are bounded on each
wavelet subspace.

\begin{example}\rm
Let us consider two symbols on $\mathbb{R}_+$ in the form $$a(v) =
v^{-\beta} \sin v^{-\alpha}, \quad \beta\in (0,1),\,\, \alpha \geq
\beta,$$ and
$$b(v) = v^{\tau} \sin v^{-\alpha}, \quad \tau\in (0,\beta).$$
According to Example~\ref{exmunbounded2}, for the unbounded symbol
$a(v)$ the Calder\'on-Toeplitz operator $T_a^{(k)}$ is bounded for
each $k\in\mathbb{Z}_+$. Since the symbol $b(v)\in C[0,+\infty]$,
then the Calder\'on-Toeplitz operator $T_b^{(k)}$ is bounded for
each $k\in\mathbb{Z}_+$ as well. Put $$ c(v) = a(v)b(v) =
\frac{v^{-\delta}}{2} - \frac{v^{-\delta}}{2} \cos 2v^{-\alpha} =
c_1(v)+c_2(v),$$ where $\delta=\beta-\tau \in (0,1)$. Clearly,
$c(v)$ is an unbounded symbol. However, $T_{c_2}^{(k)}$ is bounded
for each $k\in\mathbb{Z}_+$ (analogously as in
Example~\ref{exmunbounded2} replacing $\sin$ by $\cos$). Since
$$\theta_{c_1}(v) = \inf_{t\in (0,v)} \frac{1}{2 t^{\delta}} =
\frac{1}{2 v^{\delta}} \to +\infty,
\,\,\,\,\textrm{as}\,\,\,\,v\to 0,$$ then by Theorem~\ref{thm3}
the Calder\'on-Toeplitz operator $T_{c_1}^{(k)}$ is unbounded for
each $k\in\mathbb{Z}_+$. Thus, the Calder\'on-Toeplitz operator
$T_{ab}^{(k)}$ is unbounded on $A^{(k)}$ for each
$k\in\mathbb{Z}_+$. Moreover, this result shows that the
\textit{semi-commutator} $$\left[T_a^{(k)}, T_b^{(k)}\right) =
T_a^{(k)}T_b^{(k)} - T_{ab}^{(k)}$$ \textit{is not compact}. This
interesting feature and the algebras of Calder\'on-Toeplitz
operators will be considered elsewhere.
\end{example}


\section{Calder\'on-Toeplitz operators with unbounded symbols as uniform limits of
Calder\'on-Toeplitz operators with bounded symbols}
\label{section3}

In connection with the obtained results we now show how
Calder\'on-Toeplitz operators with unbounded symbols given in
Example~\ref{exmunbounded2} can appear as uniform limits of
Calder\'on-Toeplitz operators with bounded symbols.

\begin{theorem}
The Calder\'on-Toeplitz operator $T_a^{(k)}$ with symbol
$$a(v) = v^{-\beta} \sin v^{-\alpha}, \quad v\in\mathbb{R}_+,$$ where $\beta\in (0,1)$ and $\alpha>\beta$,
belongs to the $C^*$-algebra generated by Calder\'on-Toeplitz
operators with smooth bounded symbols on each wavelet subspace
$A^{(k)}$.
\end{theorem}

\proof Consider the sequence $\vartheta_n=(\pi n)^{-1/\alpha}$,
$n\in\mathbb{N}$, of zeros of function $a=a(v)$ and define the
sequence $$a_n(v)=\begin{cases} a(v), & v\in [\vartheta_n,
+\infty),
\\ 0, & v\in [0,\vartheta_n).
\end{cases}$$ Each symbol $a_n(v)$ is bounded and continuous.
Further each $a_n(v)$ can be uniformly approximated by smooth
symbols, and thus belongs to the $C^*$-algebra generated by
Calder\'on-Toeplitz operators with smooth bounded symbols.
According to Theorem~\ref{CTO1} the Calder\'on-Toeplitz operator
$T_a^{(k)}$ acting on $A^{(k)}$ is unitarily equivalent to the
multiplication operator $\gamma_{a,k}I$ acting on
$L_2(\mathbb{R}_+)$, where the function $\gamma_{a,k}$ is given
by~(\ref{gamma1}). Thus, {\setlength\arraycolsep{2pt}
\begin{eqnarray*}
\|T_a^{(k)} - T_{a_n}^{(k)}\| & = & \|T_{a-a_n}^{(k)}\| =
\sup_{\xi\in\mathbb{R}_+} |\gamma_{(a-a_n),k}(\xi)| \\ & = &
\sup_{\xi\in\mathbb{R}_+} \left|2\xi \int_0^{\vartheta_n}
a(v)\ell_k^2(2v\xi)\,\mathrm{d}v\right| \\ & = &
\sup_{\xi\in\mathbb{R}_+} \Biggl|2\xi
C_a^{(1)}(\vartheta_n)\ell_k^2(2\vartheta_n\xi) + 4\xi^2
\int_0^{\vartheta_n} C_a^{(1)}(v) \ell_k^2(2v\xi)\,\mathrm{d}v \\
& + & 8\xi^2 \int_0^{\vartheta_n} C_a^{(1)}(v) e^{-2v\xi}
L_k(2v\xi)L_{k-1}^{(1)}(2v\xi)\,\mathrm{d}v\Biggr|,
\end{eqnarray*}}where integration by parts has been used. Since $$C_a^{(1)}(v) = \int_0^v a(t)\,\mathrm{d}t =
\frac{v^{\alpha-\beta+1}}{\alpha}\cos v^{-\alpha} +
\mathcal{O}(v^{2\alpha-\beta+1}),\,\,\,\,\textrm{as}\,\,\,\,v\to
0,$$ see Example~\ref{exmunbounded2}, then
$$|C_a^{(1)}(v)| \leq \textrm{const}\, v^{\alpha-\beta+1},$$ where
"const" does not depend on $v\in(0,1)$, and thus the
Calder\'on-Toeplitz operator $T_a^{(k)}$ is bounded on $A^{(k)}$
for each $k\in\mathbb{Z}_+$. Then {\setlength\arraycolsep{2pt}
\begin{eqnarray*}
\|T_{a-a_n}^{(k)}\| & \leq &
\textrm{const}\,\sup_{\xi\in\mathbb{R}_+} \Biggl( 2\xi
|C_a^{(1)}(\vartheta_n)|\ell_k^2(2\vartheta_n\xi) + 4\xi^2
\int_0^{\vartheta_n} |C_a^{(1)}(v)| \ell_k^2(2v\xi)\,\mathrm{d}v
\\ & + & 8\xi^2 \int_0^{\vartheta_n} |C_a^{(1)}(v)| e^{-2v\xi}
|L_k(2v\xi)L_{k-1}^{(1)}(2v\xi)|\,\mathrm{d}v\Biggr) \\ & \leq &
\textrm{const}\,\sup_{\xi\in\mathbb{R}_+}
\vartheta_n^{\alpha-\beta}
(2\vartheta_n\xi)\ell_k^2(2\vartheta_n\xi)
\\ & + & \textrm{const}\,\sup_{\xi\in\mathbb{R}_+} 4\xi^2
\int_0^{\vartheta_n} v^{\alpha-\beta+1}
\ell_k^2(2v\xi)\,\mathrm{d}v \\ & + &
\textrm{const}\,\sup_{\xi\in\mathbb{R}_+} 8\xi^2
\int_0^{\vartheta_n}
v^{\alpha-\beta+1} e^{-2v\xi} |L_k(2v\xi)L_{k-1}^{(1)}(2v\xi)|\,\mathrm{d}v \\
& = & I_1 + I_2 + I_3.
\end{eqnarray*}}Since $\sup\limits_{\xi\in\mathbb{R}_+}
(2\vartheta_n\xi)\ell_k^2(2\vartheta_n\xi)<+\infty$, then $I_1\leq
q_1(k) \vartheta_n^{\alpha-\beta}$. To evaluate $I_2$ we use the
estimate~(\ref{upperbound2}), and we have
{\setlength\arraycolsep{2pt}
\begin{eqnarray*}
I_2 & = & \textrm{const}\,\sup_{\xi\in\mathbb{R}_+}
\int_0^{\vartheta_n} v^{\alpha-\beta} (2v\xi) \ell_k^2(2v\xi)
2\xi\,\mathrm{d}v \\ & \leq &
\textrm{const}\,\vartheta_n^{\alpha-\beta}
\sup_{\xi\in\mathbb{R}_+} \int_0^{\vartheta_n} (2v\xi)
\ell_k^2(2v\xi) 2\xi\,\mathrm{d}v \\ & = &
\textrm{const}\,\vartheta_n^{\alpha-\beta}
\sup_{\xi\in\mathbb{R}_+} \int_0^{2\vartheta_n\xi}
x\ell_k^2(x)\,\mathrm{d}x \\ & \leq &
\textrm{const}\,\vartheta_n^{\alpha-\beta} \int_{\mathbb{R}_+} \Lambda_{1,k,k}^{(0,0)}(x)\,\mathrm{d}x \\
& \leq & q_2(k)\vartheta_n^{\alpha-\beta}.
\end{eqnarray*}}Similarly for $I_3$ we get {\setlength\arraycolsep{2pt}
\begin{eqnarray*}
I_3 & = & 2\,\textrm{const}\,\sup_{\xi\in\mathbb{R}_+}
\int_0^{\vartheta_n} v^{\alpha-\beta} (2v\xi) e^{-2v\xi}
|L_k(2v\xi)L_{k-1}^{(1)}(2v\xi)| 2\xi\,\mathrm{d}v
\\ & \leq & 2\,\textrm{const}\,\vartheta_n^{\alpha-\beta}
\sup_{\xi\in\mathbb{R}_+} \int_0^{\vartheta_n} (2v\xi) e^{-2v\xi}
|L_k(2v\xi)L_{k-1}^{(1)}(2v\xi)| 2\xi\,\mathrm{d}v \\ & = &
2\,\textrm{const}\,\vartheta_n^{\alpha-\beta}
\sup_{\xi\in\mathbb{R}_+} \int_0^{2\vartheta_n\xi} xe^{-x}
|L_k(x)L_{k-1}^{(1)}(x)|\,\mathrm{d}x
\\ & \leq & 2\,\textrm{const}\,\vartheta_n^{\alpha-\beta} \int_{\mathbb{R}_+}
\Lambda_{1,k,k-1}^{(0,1)}(x)\,\mathrm{d}x \\ & \leq &
q_3(k)\vartheta_n^{\alpha-\beta}.
\end{eqnarray*}}Thus, $$\|T_a^{(k)} - T_{a_n}^{(k)}\| \leq q(k) \vartheta_n^{\alpha-\beta},$$
where the constant $q(k)$ depends on $k$, but does not depend on
$n$, and $\vartheta_n\to 0$ whenever $n\to+\infty$. \qed

It seems to be natural to ask whether the boundedness of
Calder\'on-Toeplitz operator (and by Corollary~\ref{corgamma} the
boundedness of corresponding function $\gamma_{\cdot}$) is
equivalent to the boundedness of its Wick symbol. According to the
result of Nowak~\cite{nowak3} it is true for non-negative symbols
$a$ and sufficiently smooth wavelets. Thus, we immediately have
the following result.

\begin{corollary}\rm
For a non-negative symbol $a=a(v)$ the following statements are
equivalent:
\begin{itemize}
\item[(i)] operator $T_a^{(k)}$ is bounded; \item[(ii)] the
function $\gamma_{a,k}$ is bounded; \item[(iii)] the Wick symbol
$\widetilde{a}_k$ of $T_a^{(k)}$ is bounded.
\end{itemize}
\end{corollary}

In connection with it we also mention another Nowak's result,
cf.~\cite{nowak3}, for compactness of Calder\'on-Toeplitz
operator: \textit{for a non-negative symbol the
Calder\'on-Toeplitz operator is compact if and only if its Wick
symbol tends to $0$ at infinity}. In our case of symbol $a=a(v)$
we are in a different situation because according to
Theorem~\ref{CTO1} the operator $T_a^{(k)}$ is unitarily
equivalent to a multiplication operator, and thus never compact.

\vspace{12pt} {\bf Acknowledgements:} Author has been on a
postdoctoral stay at the Departamento de Mate\-m\'aticas,
CINVESTAV del IPN (M\'exico), when writing this paper and
investigating the topics presented herein. He therefore gratefully
acknowledges the hospitality of the mathematics department of
CINVESTAV on this occasion. Author wishes especially to thank
Nikolai L. Vasilevski for having read previous versions of the
present work and for having generously shared with him his
comments and suggestions.


\vspace{5mm}

\noindent \small{Ondrej Hutn\'ik, Departamento de Matem\'aticas,
CINVESTAV del IPN, {\it Current address:} Apartado Postal 14-740,
07000, M\'exico, D.F., M\'exico
\newline {\it E-mail address:} hutnik@math.cinvestav.mx}
\newline AND \newline \noindent \small{Institute of Mathematics,
Faculty of Science, Pavol Jozef \v Saf\'arik University in Ko\v
sice, Jesenn\'a 5, 040~01 Ko\v sice, Slovakia,
\newline {\it E-mail address:} ondrej.hutnik@upjs.sk}

\end{document}